\documentclass[aps,prd,12point,nofootinbib,superscriptaddress]{revtex4} %{revtex4-1}
\def\theequation{\arabic{section}.\arabic{equation}}
\usepackage{psfrag}
\usepackage{subfigure}
\usepackage{color}
\usepackage{mathrsfs}
\usepackage{graphicx}
\usepackage{amssymb, bm}
\usepackage{amsmath, amsthm}
\usepackage{epstopdf}
\usepackage{hyperref}
\usepackage{enumerate}
\usepackage{longtable}
\usepackage{amsmath}

\usepackage{amsmath}  % needed for \tfrac, \bmatrix, etc.
\usepackage{amsfonts} % needed for bold Greek, etc.
\usepackage{graphicx}

\newcommand{\be}{\begin{equation}}
\newcommand{\ee}{\end{equation}}

\begin{document}
\def\theequation{\arabic{section}.\arabic{equation}}

\title{One-parameter Darboux-deformed Fibonacci numbers}
\author{H.C. Rosu}
\email[hcr@ipicyt.edu.mx; ORCID: 0000-0001-5909-1945]{}
\affiliation{Instituto Potosino de Investigaci\'on Cient\'{\i}fica y Tecnol\'ogica, Camino a la Presa San Jos\'e 2055,
Col. Lomas 4a Secci\'on, San Luis Potos\'{\i}, 78216 S.L.P., Mexico
}

\author{S.C. Mancas}
\email[mancass@erau.edu, ORCID: 0000-0003-1175-6869]{}
\affiliation{Department of Mathematics,
Embry-Riddle Aeronautical University, Daytona Beach, FL 32114-3900, USA}

%\date{\today}

\bigskip
\bigskip
\begin{abstract}

One-parameter Darboux deformations are effected for the simple ODE satisfied by the continuous generalizations of the Fibonacci sequence recently discussed by Faraoni and Atieh [Symmetry 13, 200 (2021)], who promoted a formal analogy with the Friedmann equation in the FLRW homogeneous cosmology.
The method allows the introduction of deformations of the continuous Fibonacci sequences, hence of Darboux-deformed Fibonacci (non integer) numbers.
Considering the same ODE as a parametric oscillator equation, the Ermakov-Lewis invariants for these sequences are also discussed.
\end{abstract}

%\pacs{}
%\keywords{}

\maketitle

\section{Introduction}
\setcounter{equation}{0}
\label{sec:1}

If the famous Binet formula for the Fibonacci numbers is written using exponentials
\be\label{binetd}
F_n=\frac{e^{\tilde{\varphi}n}-(-1)^ne^{-\tilde{\varphi}n}}{\sqrt{5}}
\ee
where $\tilde{\varphi}\approx 0.4812$ is the natural logarithm of the golden ratio $\varphi \equiv \frac{1 + \sqrt{5}}{2} \approx 1.6180$, then one may think that the natural continuous generalization is the function
\be\label{binetc}
F_c(x)=\frac{e^{\tilde{\varphi}x}-\cos(\pi x)e^{-\tilde{\varphi}x}}{\sqrt{5}}~.
\ee

\begin{figure}
\begin{center}
\includegraphics[width=8.5cm]{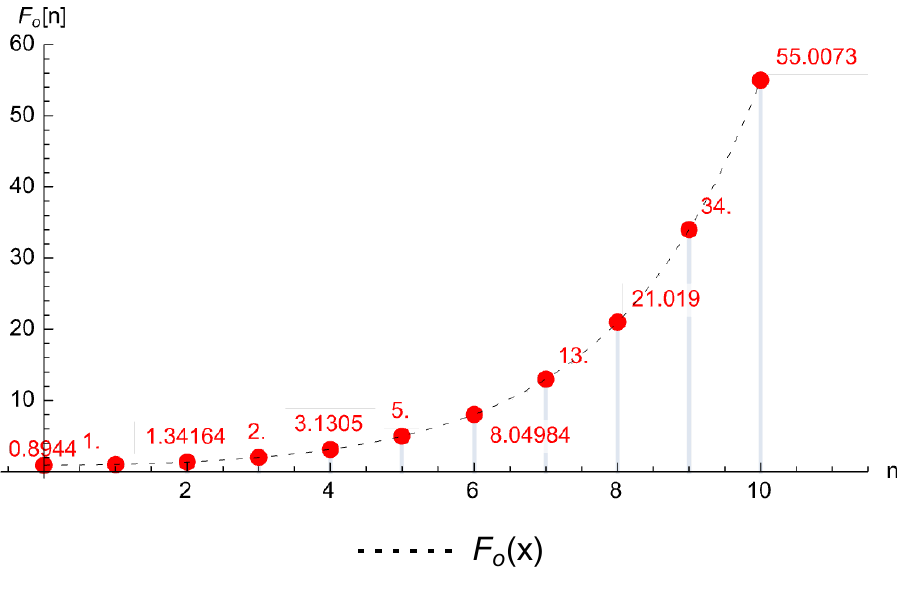}
\includegraphics[width=8.5cm]{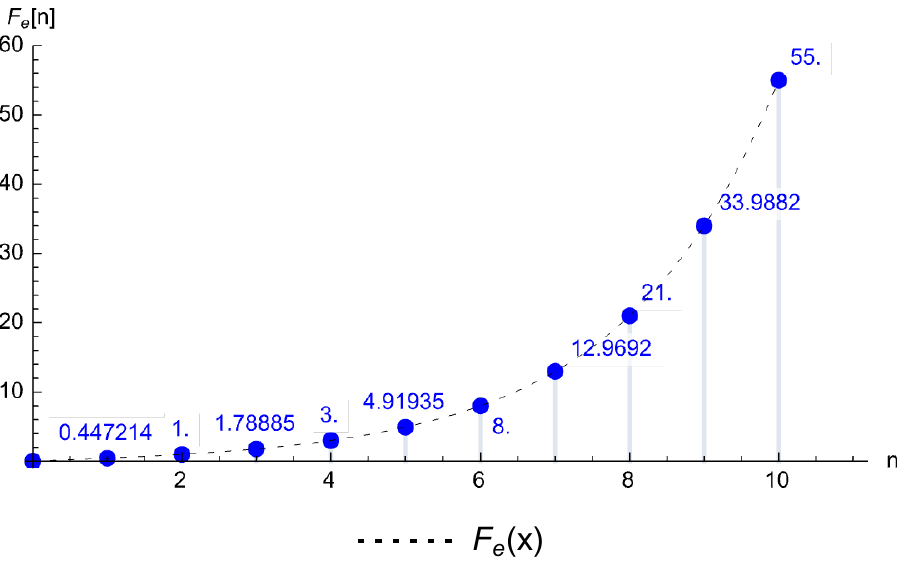}
\caption{The red  balls ($n$ is odd), and blue balls ($n$ is even) correspond to the first ten Fibonacci numbers given by the sequence \eqref{binetd}, while the dotted curves are the corresponding hyperbolic functions \eqref{bit2}.   }\label{f1}
\end{center}
\end{figure}

In a recent short note, Faraoni and Atieh \cite{FA} pointed out that the even and odd components of (\ref{binetc})
\begin{eqnarray}\label{bit2}
&F_{e}(x)=\frac{2}{\sqrt{5}}\, \sinh \left(\tilde{\varphi} x  \right)~, \\ %\nonumber\\
&F_{o}(x)=\frac{2}{\sqrt{5}}\, \cosh \left(\tilde{\varphi}  x  \right)~,   \nonumber
\end{eqnarray}
which provide Fibonacci numbers only for $x=n\in \mathbb{N}$ even and odd, respectively, see FIG.~\ref{f1},
can be treated as equivalents of the scale factors of spatially homogeneous and isotropic cosmologies.
This is not so surprising if one thinks of the strongly expanding feature of the gaps between Fibonacci's numbers generated by their recurrence relationship.
Furthermore, the fact that the even and odd Binet components are hyperbolic sine and cosine functions, respectively, reveal them as linear independent solutions to the simple hyperbolic `oscillator'
\be\label{odeF}
F_i''-\tilde{\varphi}^2F_i=0 ~,
\ee
where the subscript $i$ stands for $e$ or $o$. On the other hand, the continuous generalization (\ref{binetc}), which reproduces all of the Fibonacci numbers when $x\in \mathbb{N}$, is the solution of a similar, but nonhomogeneous differential equation.

As commented in \cite{FA}, in classical mechanics the hyperbolic oscillator can be used only
as the simplest example of an unstable mechanical system of a particle of position $F$ in the inverted harmonic oscillator potential
$V(F)=-KF^2/2$ with the elastic constant $K=\tilde{\varphi}^2$. This also shows why the cosmological analogy may be considered as more interesting.

\medskip

In addition, Faraoni and Atieh obtained Lagrangian and Hamiltonian formulations
for the Fibonacci ODEs (\ref{odeF}) and further derived from these an invariant of the
(discrete) Fibonacci sequence.

\medskip

Here, we construct the Darboux deformed counterpart of equation (\ref{odeF}) based on the general Riccati solution which introduces a one-parameter distortion of the Fibonacci numbers, especially of the small ones for small values of the parameter and becomes exponentially vanishing when advancing along the sequence. The parametric Darboux scheme has long been used in supersymmetric quantum mechanics ever since it has been introduced in that context thirty-five years ago \cite{Mielnik,Fernandez}; for a more recent paper, see \cite{RMC}, and for a close paper to the present one, see \cite{rr98}.

We also use the Ermakov nonlinear equation associated to (\ref{odeF}) to introduce the Ermakov invariant for these continuous Fibonacci sequences.

\bigskip

\section{The parametric Darboux deformation of $F_{e,o}(x)$}

In this section, we first present the parametric Darboux deformation in a sufficiently general way to provide its key points, and next move to the two Fibonacci cases.
The operatorial form of (\ref{odeF})
\be\label{SchrFib}
\left(D^2- f(x) \right) F=0~,\qquad D=\frac{d^2}{dx^2}~,
\ee
can be factored in the form
\be\label{SchrFib1}
(D-\Phi)(D+\Phi)F\equiv \big[D^2-(-\Phi'+\Phi^2)\big]F=0~,
\ee
with the negative logarithmic derivative, $\Phi=-F'/F$ of a solution $F$ of (\ref{SchrFib}). %which are
In the $\Phi$ variable, (\ref{SchrFib}) is turned into a Riccati equation, which in mathematical terminology is the well-known reduction of order of linear second-order ODEs,
\be\label{RicFib}
-\Phi'+\Phi^2=f(x)
\ee
and if one knows a Riccati solution $\Phi$, one can proceed viceversa to obtain the solution of the corresponding linear second-order ODE from $F={\rm exp}(-\int^x \Phi)$.

\medskip

Suppose now that like in supersymmetric quantum mechanics, we revert the factoring to obtain the partner equation
\be\label{SchrFibDT}
(D+\Phi)(D-\Phi)\hat{F}\equiv \big[D^2-(\Phi'+\Phi^2)\big]\hat{F}=0~,
\ee
which in mathematical terms is known as the Darboux-transformed equation of (\ref{SchrFib1}). The generic interesting fact of the reverted factorization which leads to the Darboux partner equation is that it is not unique \cite{Mielnik,Fernandez,RMC}. Indeed, by adding a function $1/u(x)$ to the logarithmic derivative in the factorization brackets, we obtain
\be\label{SchrFibDTp}%2.8
\left(D+\Phi+\frac{1}{u}\right)\left(D-\Phi-\frac{1}{u}\right)\tilde{F}\equiv \big[D^2-(\Phi'+\Phi^2)\big]\tilde{F}
-\big[\left(-u'+2\Phi u+1\right)/u^2\big]\tilde{F}=0~.
\ee
The latter formula shows that if the function $u$ satisfies the first-order differential equation
\be\label{SchrFibDTl}
-u'+2\Phi u+1=0
\ee
then the same Darboux-tranformed equation is obtained. Thus, $\Phi_{g}=\Phi+1/u$ cannot be but the general solution in the form of the Bernoulli ansatz for the Riccati equation corresponding to the Darboux-transformed equation since in this ansatz the linear equation (\ref{SchrFibDTl}) is the basic equation fulfilled by Bernoulli's function $u$. Moreover, the following regrouping
\be\label{SchrFibDTp1}
\bigg[D^2-\left(\Phi'-\frac{u'}{u^2}\right)-\left(\tilde{\Phi}^2+2\frac{\Phi}{u}+\frac{1}{u^2}\right)\bigg]\tilde{F}=0
\ee
can be written as
\be\label{SchrFibDTp2}
\big[D^2-\left(\Phi'_{g}+\Phi^2_{g}\right)\big]\tilde{F}=0~,
\ee
which seems to be not relevant unless one thinks about the partner equation
\be\label{SchrFibDTp3}
\big[D^2-\left(-\Phi'_{g}+\Phi^2_{g}\right)\big]{\cal F}=0~.
\ee
In the latter case, the explicit calculation leads to
\be\label{SchrFibDTp4}
\big[D^2-\left( f(x)+4\Phi u^{-1}+2u^{-2}\right)\big]{\cal F}=0~,
\ee
which represents a one-parameter family of equations that have the same Darboux-transformed partner, the running parameter of the family being the integration constant that occurs in the
Bernoulli's function $1/u$. The latter can be easily obtained by the integration of (\ref{SchrFibDTl}) in the form
%............................
\be\label{SchrFibDTp5}%2-15
\frac{1}{u}=\frac{e^{-2\int^x \Phi(x')dx'}}{\gamma +\int^x e^{-2\int^x \Phi(x')dx'}dx}=\frac{d}{dx}\ln\left(\gamma +\int^x e^{-2\int^x \Phi(x')dx'}\right)~.
\ee
The solution ${\cal F}$ is a parametric Darboux-deformed partner of $F$.
Its explicit form in terms of $F$ can be obtained from
%...........................
\be\label{SchrFibDTp6}%2-15
{\cal F}=e^{-\int^x \Phi_{g} ds}=e^{-\int^x \Phi ds}e^{-\int^x\frac{d}{dx}\ln\left(\gamma +\int^x e^{-2\int^x \Phi(s)ds}\right)}
=\frac{F}{\gamma +\int^x F^2(s)ds}~.
\ee

\medskip

We are now ready to apply this simple mathematical scheme to the two continuous Fibonacci sequences for which $f(x)=\tilde{\varphi}^2$. We start
with the odd case followed by the even case.

\medskip

\subsection{The odd deformed case}
In the odd continuous Fibonacci case, $\Phi_o=-\tilde{\varphi}\tanh(\tilde{\varphi}x)$ and the Darboux pair of
Riccati equations have the form
%.............................
\be\label{ric2}%2.17
-\frac{d\Phi_o}{dx}+\Phi_o^2=\tilde{\varphi}^2~, \qquad \frac{d\Phi_o}{dx}+\Phi_o^2=\tilde{\varphi}^2[1-2{\rm sech}^2(\tilde{\varphi}x)]~,
\ee
%.............................
where the soliton-type free term in the right hand side of the second Riccati equation is the non-parametric Darboux deformation of $\tilde{\varphi}^2$.
For this case, the parametric Darboux partner equation (\ref{SchrFibDTp4}) takes the form
%.........................
\be\label{peqodd}%2.17
\left(\frac{d^2}{dx^2}- \tilde{\varphi}_{\gamma,o}^2\right) {\cal F}_o=0~, \qquad
\tilde{\varphi}_{\gamma,o}^2=\tilde{\varphi}^2\left(1-\frac{8\sinh(2\tilde{\varphi}x)}{2\tilde{\varphi}(2\gamma+x)+\sinh(2\tilde{\varphi}x)}+
\frac{32\cosh^4(\tilde{\varphi}x)}{(2\tilde{\varphi}(2\gamma+x)+\sinh(2\tilde{\varphi}x))^2}\right)
\ee
%........................
with two independent solutions given by
\begin{align}%2.18
&{\cal F}_o(\tilde{\varphi}x)=A_oF_o~, \quad  %=\frac{\cosh(\tilde{\varphi}x)}{\gamma+\frac{x}{2}+\frac{\sinh(2\tilde{\varphi}x)}{4\tilde{\varphi}}}=
A_o=\frac{2\sqrt{5}\tilde \varphi}{2 \tilde \varphi(2 \gamma +x)+\sinh (2\tilde  \varphi  x)}~, \label{psolA}\\ %F_o~,
&{\cal G}_o(\tilde{\varphi}x)=B_oF_o~,  \quad
%\frac{W \cosh (\tilde \varphi  x) \left[4 \tilde \varphi ^2 (2 \gamma +x)^2 \tanh (\tilde \varphi  x)-2 \tilde \varphi  x+\sinh (2 \tilde \varphi  x)\right]}{4\tilde  \varphi ^2 \left[2 \tilde \varphi  (2 \gamma +x)+\sinh (2\tilde  \varphi  x)\right]}~\\
B_o=\frac{\sqrt{5}W\left[4 \tilde \varphi ^2 (2 \gamma +x)^2 \tanh (\tilde \varphi  x)-2 \tilde \varphi  x+\sinh (2 \tilde \varphi  x)\right]}{8\tilde  \varphi ^2 \left[2 \tilde \varphi  (2 \gamma +x)+\sinh (2\tilde  \varphi  x)\right]}~.\label{psolB} %F_o~.
\end{align}
The plots of the factors $A_o$ and $B_o$ in FIG.~2 shows that ${\cal F}_o$ cannot serve as a deformed Fibonacci sequence since $A_o$ goes rapidly to zero.
Instead, $B_o$ goes to a small constant value and as such it turns ${\cal G}_o(\tilde{\varphi}x)$ into the appropriate deformed odd counterpart that follows approximately the odd Fibonacci sequence. We present
plots of the Darboux-deformed ${\cal G}_o(\tilde{\varphi}x)$ sequence for three values of the deformation parameter $\gamma$ in FIG.~\ref{f2}.
The effect of the deformation is stronger at the beginning of the sequence, but quickly vanishes after the first few Fibonacci numbers.

\begin{figure}
\begin{center}
\includegraphics[width=8.5cm]{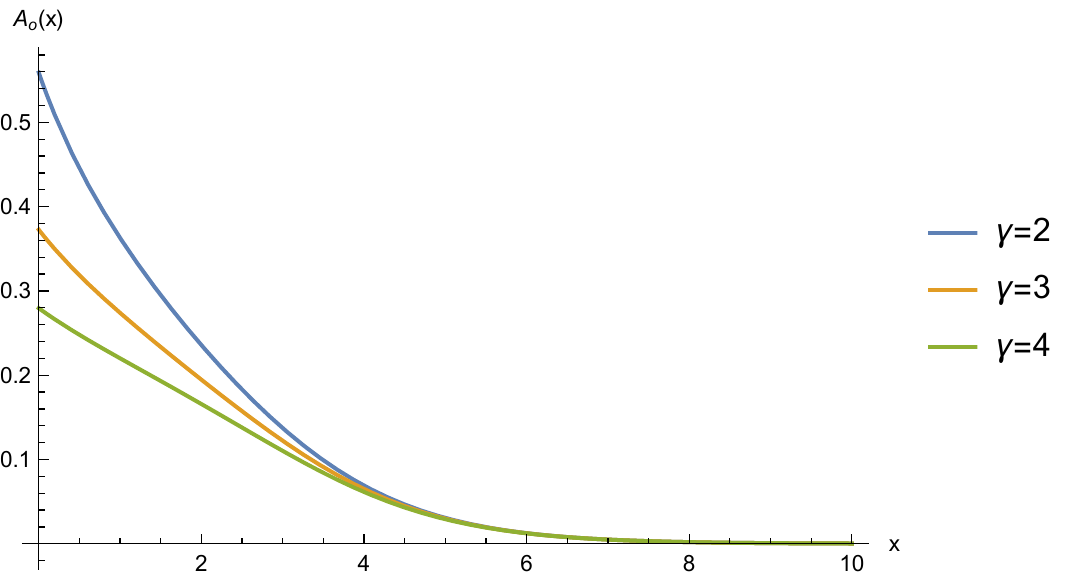}
\includegraphics[width=8.5cm]{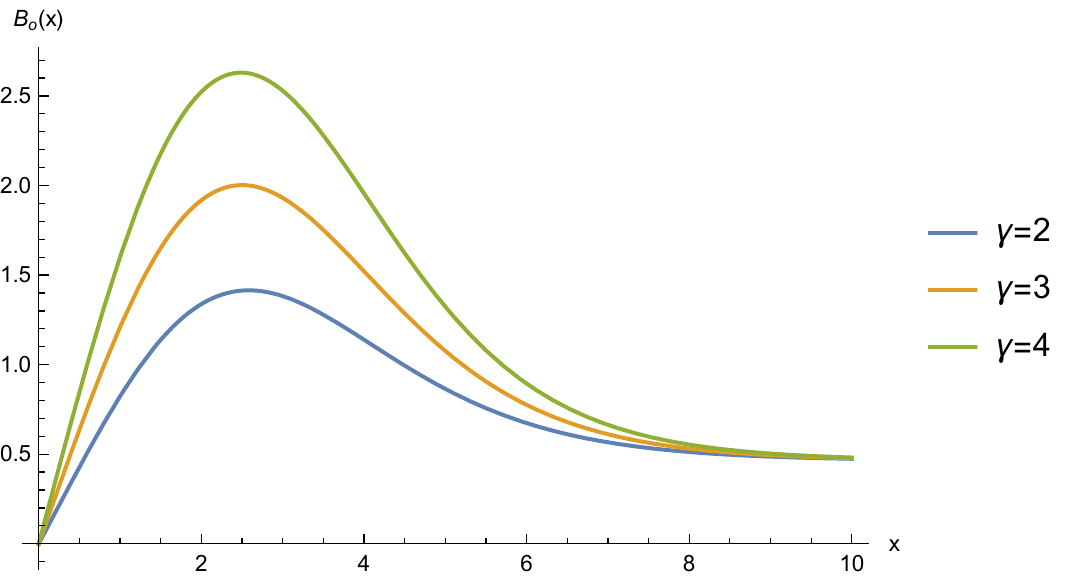}
\caption{The factors $A_o$ and $B_o$ vs $x$.   }\label{f1}
\end{center}
\end{figure}

\begin{figure}
\begin{center}
\includegraphics[width=16cm]{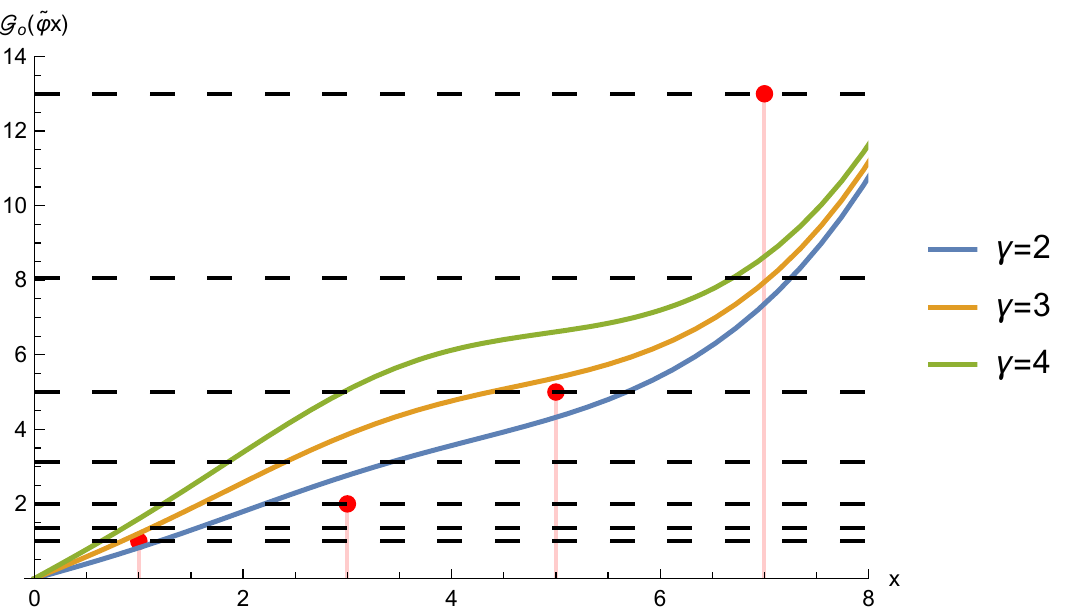}
\caption{The red balls correspond to the first four Fibonacci numbers in the odd Fibonacci sequence $F_o$. The curves are the corresponding Darboux-deformed sequences ${\cal G}_o$ for three values of the deformation parameter.}\label{f2}
\end{center}
\end{figure}

The Darboux-deformed Fibonacci `odd' numbers that can be introduced are real numbers located at
\begin{align}%\label{ddodd}
%&\Delta_{A_o}(\tilde{\varphi} x)=\Bigg[\frac{2\sqrt{5}\tilde \varphi}{2 \tilde \varphi(2 \gamma +x)+\sinh (2\tilde  \varphi  x)}-1\Bigg]
%F_o(\tilde{\varphi}x) \label{DAo}\\
&\Delta_{B_o}(\tilde{\varphi} x)=\Bigg[\frac{\sqrt{5}W\left[4 \tilde \varphi ^2 (2 \gamma +x)^2 \tanh (\tilde \varphi  x)-2 \tilde \varphi  x+\sinh (2 \tilde \varphi  x)\right]}{8\tilde  \varphi ^2 \left[2 \tilde \varphi  (2 \gamma +x)+\sinh (2\tilde  \varphi  x)\right]}-1\Bigg]F_o(\tilde{\varphi}x) \label{DBo}
\end{align}
with respect to the corresponding odd Fibonacci numbers at $x=2m-1$, $m=1,2,...$. In Table 1, we present these Darboux shifts for the first four odd Fibonacci numbers.

\begin{center}
\begin{tabular}{|c|c|c|c|}

  \hline $x=2m-1$ & $\Delta_{B_o}$ $(\gamma=2)$ & $\Delta_{B_o}$ $(\gamma=3)$ & $\Delta_{B_o}$ $(\gamma=4)$\\

 \hline F$_1$ = 1 & -0.17634 & 0.209938 & 0.601617\\
 \hline F$_3$ = 2& 0.764837  & 1.86197& 3.05822 \\
 \hline F$_5$ = 5& -0.678174 & 0.381427  & 1.6148\\
 \hline F$_7=$13 & -5.63824 & -5.05354 & -4.37156\\
  \hline
\end{tabular}
\hspace{3.5cm}
\end{center}
%\hspace*{3.5cm}
{Table 1. The first four Fibonacci numbers $\{1,2,5,13\}$ in the odd sequence and their Darboux shifts.}

\medskip

\subsection{The even deformed case}
\begin{figure}[h!]
\begin{center}
\includegraphics[width=16cm]{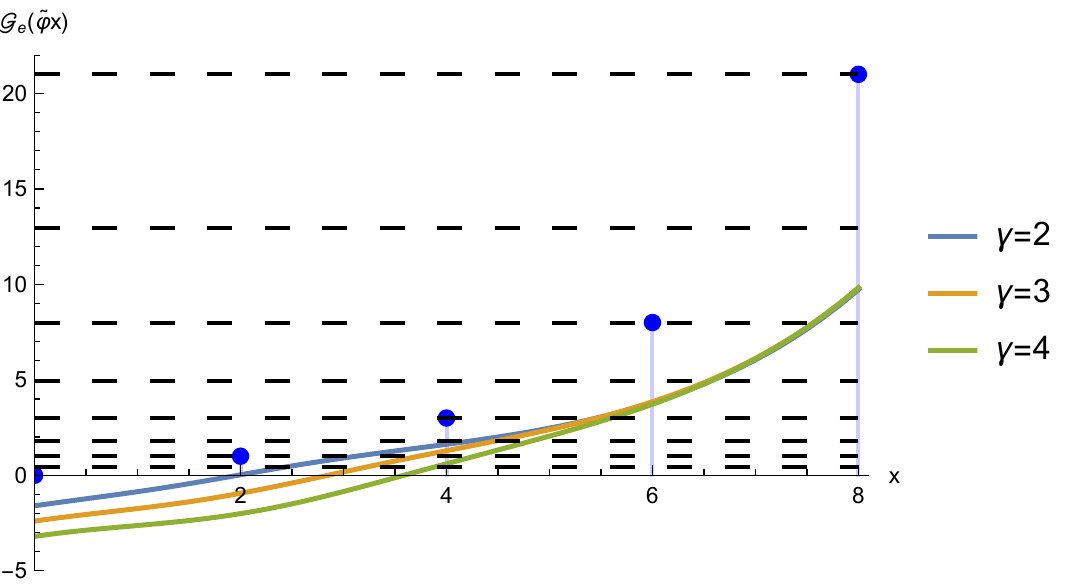}
\caption{ The blue balls correspond to the first five Fibonacci numbers in the even Fibonacci sequence $F_e$. The curves are the corresponding Darboux-deformed sequences ${\cal G}_e$ for the same three values of the deformation parameter as in the previous figure.}\label{f3}
\end{center}
\end{figure}

In the even continuous Fibonacci case, $\Phi_e=-\tilde{\varphi}\coth(\tilde{\varphi}x)$ and the Darboux pair of
Riccati equations have the form
\be\label{ric2}
-\frac{d\Phi_e}{dx}+\Phi_e^2=\tilde{\varphi}^2~, \qquad \frac{d\Phi_e}{dx}+\Phi_e^2=\tilde{\varphi}^2[1+2{\rm cosech}^2(\tilde{\varphi}x)]~.
\ee
The parametric Darboux partner equation reads
\be\label{peqeven}
\left(\frac{d^2}{dx^2}- \tilde{\varphi}_{\gamma,e}^2\right) {\cal F}_e=0~, \qquad
\tilde{\varphi}_{\gamma,e}^2=\tilde{\varphi}^2\left(1-\frac{8\sinh(2\tilde{\varphi}x)}{2\tilde{\varphi}(2\gamma-x)+\sinh(2\tilde{\varphi}x)}+
\frac{32\sinh^4(\tilde{\varphi}x)}{(2\tilde{\varphi}(2\gamma-x)+\sinh(2\tilde{\varphi}x))^2}\right)
\ee
with the two independent solutions given by
\begin{align}
&{\cal F}_e(\tilde{\varphi}x)=A_eF_e~, \quad A_e=\frac{2\sqrt{5}\tilde \varphi}{2 \tilde \varphi(2 \gamma -x)+\sinh (2\tilde  \varphi  x)}~,\label{pseA}\\
&{\cal G}_e(\tilde{\varphi}x)=B_eF_e~, \quad
B_e=\frac{\sqrt{5}W \left[-4 \tilde\varphi ^2 (x-2 \gamma )^2 \coth (\tilde\varphi  x)+2 \tilde\varphi  x+\sinh (2\tilde \varphi  x)\right]}{8 \tilde\varphi ^2 \left[2\tilde\varphi(2\gamma -x)+\sinh (2 \tilde\varphi  x)\right]} \label{pseB}~.
\end{align}
Again, we find that only the ${\cal G}_e(\tilde{\varphi}x)$ solution follows relatively close the even Fibonacci sequence. Thus, we display plots of it for the same three values of the Darboux parameter as in the odd case in FIG.~\ref{f3}.
Similarly to the odd case, we can introduce Darboux-deformed Fibonacci `even' numbers according to %(\ref{pseA}) and
(\ref{pseB}), which are real numbers located at
\begin{align}%\label{ddodd}
%&\Delta_{A_e}(\tilde{\varphi} x)=\Bigg[\frac{2\sqrt{5}\tilde \varphi}{2 \tilde \varphi(2 \gamma -x)+\sinh (2\tilde  \varphi  x)}-1
%\Bigg]F_e(\tilde{\varphi}x) \label{DAo}\\
&\Delta_{B_e}(\tilde{\varphi} x)=\Bigg[\frac{\sqrt{5}W \left[-4 \tilde\varphi ^2 (x-2 \gamma )^2 \coth (\tilde\varphi  x)+2 \tilde\varphi  x+\sinh (2\tilde \varphi  x)\right]}{8 \tilde\varphi ^2 \left[2\tilde\varphi(2\gamma -x)+\sinh (2 \tilde\varphi  x)\right]}-1
\Bigg]F_e(\tilde{\varphi}x) \label{DBo}
\end{align}
with respect to the corresponding even Fibonacci number at $x=2m$, $m=0,1,...$. In Table 2, we present these Darboux shifts for the first five Fibonacci numbers in the even positions 0,2,4,6, and 8.

\begin{center}
\begin{tabular}{|c|c|c|c|}

  \hline $x=2m$ & $\Delta_{B_e}$ $(\gamma=2)$ & $\Delta_{B_e}$ $(\gamma=3)$ & $\Delta_{B_e}$ $(\gamma=4)$\\

 \hline F$_0$ = 0&  -1.600 &  -2.400 & -3.200 \\
 \hline F$_2$ = 1&  -0.972878 &  -1.94204 & -3.00855 \\
 \hline F$_4$ = 3 & -1.37741 &  -1.71251 & -2.3984\\
 \hline F$_6$ = 8 & -4.1898  & -4.14927 & -4.27983\\
 \hline F$_8$ = 21& -11.271 & -11.1894  & -11.1738\\
  \hline
\end{tabular}
\hspace{3.5cm}
\end{center}
%\hspace*{3.5cm}
Table 2. The first five Fibonacci numbers $\{0,1,3,8,21\}$ in the even sequence and their Darboux shifts.

\bigskip

\section{The Ermakov-Lewis invariant}\label{S2}

In this section, we consider (\ref{odeF}) as a parametric oscillator equation in the particular cases of constant coefficients, which has been studied in more detail in \cite{maro}. For the $x$-axis being a number axis, there is no problem to consider it as a temporal axis.
It is a well-known result that one can use two given linear independent solutions, $u_1$, and $u_2$, of the parametric oscillator equation,
%........3
\begin{equation}\label{a9}
u'' +\omega^2(x)u=0,
\end{equation}
to build a particular solution of the corresponding nonlinear Ermakov-Pinney equation
%........4
\begin{equation}\label{a10}
v''+\omega^2(x)v+kv^{-3}=0
\end{equation}
by means of Pinney's formula \cite{P}
%.........5
\begin{equation}\label{a11}
v(x)=\sqrt{u_1^2-\frac{ku_2^2}{W^2}}~,
\end{equation}
where $W$ is the Wronskian of the two solutions $u_1$, and $u_2$.

\medskip

The Fibonacci case corresponds to $\omega^2(x)=-\tilde{\varphi}^2$ with the Ermakov-Pinney equation
%........7
\begin{equation}\label{a13}
v''-\tilde{\varphi}^2v+kv^{-3}=0~.
\end{equation}
Such an equation has been called the simplest EP (SEP) equation in \cite{maro}.
From (\ref{a11}), by taking $u_1=F_{e}(x), ~u_2=F_{o}(x),$ the particular solution with $W=-\frac{4 \tilde \varphi}{5} $ can be written immediately in the form \cite{maro}
%........9
\be\label{a14v}
v(x)=\sqrt{\frac{4}{5} \sinh ^2(\tilde\varphi  x)-\frac{5 k \cosh ^2(\tilde\varphi  x)}{4 \tilde\varphi ^2}}~.
\ee

In general, the Ermakov invariant is defined as
\be\label{einv}
{\cal I}=\frac{1}{2}\bigg[-k\Big(\frac{u}{v}\Big)^2+(u'v-uv')^2\bigg]~,
\ee and  by taking the linear superposition $u=M u_1+N u_2$ \cite{mr2014}, with $v$ from \eqref{a14v}, it becomes
\be\label{einv}
{\cal I}=\frac{16 N ^2 \tilde \varphi ^2}{25}-M ^2 k \equiv const.
\ee

\medskip

For the parametric Darboux-deformed counterparts, the Ermakov-Pinney equations are

\begin{equation}\label{vcal}
{\cal V}_{o}''+\tilde{\varphi}_{\gamma,o}^2{\cal V}_{o}+k_1{\cal V}_{o}^{-3}=0~, \quad
{\cal V}_{e}''+\tilde{\varphi}_{\gamma,e}^2{\cal V}_{e}+k_2{\cal V}_{e}^{-3}=0~.
\end{equation}

In the odd case, we take  $u_1={\cal F}_o(\tilde{\varphi}x)$, and $u_2={\cal G}_o(\tilde{\varphi}x)$ which also satisfy \eqref{peqodd}. Using \eqref{a11}, we have
\be\label{vodd}
{\cal V}_o=\frac{ \sqrt{256 \tilde\varphi ^6-k_1 \left(4 \tilde\varphi ^2 (2 \gamma +x)^2 \tanh (\tilde\varphi  x)-2 \tilde\varphi  x+\sinh (2 \tilde\varphi  x)\right)^2}}{4 \tilde\varphi ^2 (4 \gamma \tilde \varphi +2 \tilde\varphi  x+\sinh (2 \tilde\varphi  x))}\cosh (\tilde\varphi  x)~.
\ee

In the even case, we take  $u_1={\cal F}_e(\tilde{\varphi}x)$, and $u_2={\cal G}_e(\tilde{\varphi}x)$ which satisfy \eqref{peqeven}. Using again Pinney's formula, %\eqref{a11}
we have
\be\label{veven}
{\cal V}_e=\frac{ \sqrt{256 \tilde\varphi ^6-k_2 \left(-4 \tilde\varphi ^2 (x-2 \gamma )^2 \coth (\tilde\varphi  x)+2 \tilde\varphi  x+\sinh (2 \tilde\varphi  x)\right)^2}}{4 \tilde\varphi ^2 (4 \gamma  \tilde\varphi -2 \tilde\varphi  x+\sinh (2 \tilde\varphi  x))}\sinh (\tilde\varphi  x)~.
\ee

As before, letting  $u=\tilde{M} u_1+\tilde{N} u_2$ with ${\cal V}$ from either \eqref{vodd}, or \eqref{veven}, the invariants become
\be\label{einv2}
{\cal I}_o= \tilde{N} ^2 W^2 -{\tilde M} ^2 k_1 \equiv const~, \quad {\cal I}_e= \tilde{N} ^2 W^2 -{\tilde M} ^2 k_2 \equiv const~.
\ee
This shows that the Ermakov-Lewis invariants do not depend on the Darboux deformation parameter. Of course, the superposition constants, the Wronskian, and the nonlinear coupling constants can be suitably chosen to have identical invariants.

\medskip

\section{Conclusions}
\setcounter{equation}{0}
\label{sec:5}

We have introduced one-parameter Darboux-deformed families of the continuous Fibonacci sequences.
We also discussed the Ermakov-Lewis invariant for both the continuous Fibonacci sequences and the Darboux-deformed ones.

%\newpage

\end{document}